\documentclass[11pt,a4paper,reqno]{article}
\usepackage{template}

\renewcommand{\P}{\mathbb{P}}
\newcommand{\E}{\mathbb{E}}
\newcommand{\R}{\mathbb{R}}
\newcommand{\N}{\mathbb{N}}

\newcommand{\sP}{\mathscr{P}}

\DeclareMathOperator{\cov}{Cov}

  \newcounter{iconst}

\begin{document}

\title{Dynamical noise sensitivity for the voter model}

\date{\today}
\author{Gideon Amir
  \thanks{Email: \ \texttt{gidi.amir@gmail.com}; \ Department of Mathematics, Bar-Ilan University, 5290002, Ramat Gan, Israel.}
  \and
  Omer Angel
  \thanks{Email: \ \texttt{angel@math.ubc.ca}; \ Department of Mathematics, University of British Columbia, Vancouver, BC, V6T 1Z2, Canada.}
  \and
  Rangel Baldasso
  \thanks{Email: \ \texttt{r.baldasso@math.leidenuniv.nl}; \ Mathematical Institute, Leiden University, P.O. Box 9512, 2300 RA Leiden, The Netherlands.}
  \and
  Ron Peretz
  \thanks{Email: \ \texttt{ron.peretz@biu.ac.il}; \ Economics Department, Bar-Ilan University, 5290002, Ramat Gan, Israel.}
}
\maketitle

\begin{abstract}
We study how the consensus opinion of the voter model on finite graphs varies in light of noise sensitivity with respect to both the initial opinions and the dynamics. We first prove that the final opinion is stable with respect to small perturbations of the initial configuration. Different effects are observed when a perturbation is introduced in the dynamics governing the evolution of the process, and the final opinion is noise sensitive in this case. Our proofs rely on the relationship between the voter model and coalescing random walks, and on a precise description of this evolution when we have coupled dynamics.
\end{abstract}

\section{Introduction}

In this note we consider the voter model on finite graphs.
We study the consensus opinion in light of noise sensitivity with respect to the initial condition and the dynamics.
The notion of noise sensitivity for Boolean functions was introduced in the seminal paper~\cite{bks} by Benjamini, Kalai, and Schramm, and regards the behaviour of Boolean functions when a small fraction of its entries is subject to an independent noise.
More precisely, let $\P_{p}$ denote the probability measure in $\{0,1\}^{n}$ with independent marginals, each Bernoulli$(p)$ for $p \in [0,1]$.
Given $\omega \in \{0,1\}^{n}$ distributed according to $\P_{p}$ (The Bernoulli($p$) product measure), define the $\epsilon$-perturbed configuration $\omega^{\epsilon}$ by resampling each bit of $\omega$ independently with probability $\epsilon$.
A sequence of Boolean functions $f_{n}: \{0,1\}^{n} \to \{0,1\}$ is said to be \textbf{noise sensitive} (at level $p$) if for every $\epsilon>0$ we have that $f_{n}(\omega)$ and $f_{n}(\omega^{\epsilon})$ are asymptotically independent.
Formally,
\begin{equation}
  \lim_{n} \E_{p}\big[ f_{n}(\omega)f_{n}(\omega^{\epsilon}) \big] - \E_{p}\big[ f_{n}(\omega) \big]^{2} = 0.
\end{equation}
In contrast, the sequence $(f_{n})_{n \in \N}$ is said to be \textbf{noise stable} (at level $p$) if
\begin{equation}
  \lim_{\epsilon \to 0} \sup_{n} \P_{p}\big[ f_{n}(\omega) \neq f_{n}(\omega^{\epsilon}) \big] =0,
\end{equation}
so that small noises are not enough to provoke a macroscopic change in the value of the functions $f_{n}$.

The voter model is a classical model of interacting particles \cite{liggett}.
  Consider a connected graph with $n$ vertices $G_{n}= \big(V_n, E_{n} \big)$.
  At each time $t$, each vertex $v$ has an opinion $\eta_t(v) \in S$ for some space $S$ of possible opinions.
  The case $S=\{0,1\}$ is the most commonly used case.
  At time $t=0$, each vertex of $G_{n}$ has an opinion $\eta_0(v) \in S$.
  Each edge $\{u,v\} \in E_{n}$ is then endowed with two independent Poisson processes of rate one, with one process for each orientation of the edge. All Poisson processes for the different edges are taken to be independent from each other as well as from the initial configuration $\eta_0$.
  These processes control changes in opinions in $u$ and $v$.
  Whenever a clock tick happens in the process associated with the oriented edge $(u,v)$, the vertex $v$ copies the opinion of $u$.
(See below for a more formal definition.)
One construction of the voter model, referred to as a \textbf{graphical construction} is in terms of a collection of Poisson processes associated with the edges of the graph.
This is described in detail in Section~\ref{sec:prelim}.
We denote by $\sP$ the collection of Poisson point processes used in the graphical construction of the voter model, and write $\eta_{t}(x)$ for the resulting opinion of vertex $x$ at time $t$.

  The configuration $\eta_t\in S^{V_n}$ is a Markov chain which is almost surely eventually absorbed in a constant configuration.
  The eventual constant opinion is called the \textbf{consensus opinion}.
  We denote this final opinion by $f_{n}(\eta_{0}, \sP)$,
  depending implicitly on $G_{n}$ with initial opinion $\eta_{0}$ and graphical construction given by the clocks $\sP$.
  Given $p\in(0,1)$, we say that the consensus opinion $f_{n}$ is noise stable with respect to the initial condition if when the initial opinions are sampled from $\P_p$ we have
  \begin{equation}
    \lim_{\epsilon \to 0} \sup_{n} \P_{p}\big[ f_{n}(\eta_{0}, \sP) \neq f_{n}(\eta_{0}^{\epsilon}, \sP) \big] =0.
  \end{equation}

  To study the noise sensitivity of the consensus opinion with respect to the clock ticks of the graphical construction, it is necessary to first define what we mean by an $\epsilon$-perturbation of $\sP$.
  Given $\epsilon \in [0,1]$, $\sP^{\epsilon}$ is constructed from $\sP$ in two steps: one first obtains $\sP'$ as an $\epsilon$-thinning of $\sP$. Then define
  \begin{equation}
    \sP^{\epsilon} = \sP' \cup \sP(\epsilon),
  \end{equation}
  where $\sP(\epsilon)$ is an independent collection of clock ticks with intensity $\epsilon$.
  We say that the consensus opinion is noise sensitive with respect to the clock ticks if
  \begin{equation}
    \lim_{n} \E_{p}\big[ f_{n}(\eta_{0}, \sP)f_{n}(\eta_{0}, \sP^{\epsilon}) \big] - \E_{p}\big[ f_{n}(\eta_{0}, \sP) \big]^{2} = 0.
  \end{equation}

The main result of this note is the following theorem.
\begin{theorem}\label{t:consensus_opinion}
  Let $G_n$ be any sequence of connected graphs, where $G_n$ has $n$ vertices.
  The consensus opinion in the voter model is noise stable with respect to the initial condition for all $p \in (0,1)$ and noise sensitive with respect to the clock ticks.
  Moreover, the limit
  \begin{equation}
    \limsup_{n}
    \cov \Big( f_{n}(\eta_{0}, \sP) ,f_{n}(\eta_{0}, \sP^{\epsilon_{n}}) \Big)
  \end{equation}
  equals zero if $n\epsilon_{n}\to\infty$ and is positive if $n\epsilon_{n}$ is bounded.
\end{theorem}

One way to interpret the above theorem, suggested to us by Itai Benjamini, is by viewing the voter model dynamics as a random Boolean function of the initial configuration.
The dynamics determines a (random) dictator, whose opinion is revealed when checking the initial opinions.
Our theorem then says that the choice of the dictator is sensitive with respect to noising the dynamics.
Clearly, once the dictator's identity is changed all correlation is lost.
On the other hand, given a choice of dictator, her opinion is clearly stable with respect to small perturbations.
The stability statement in the theorem is an immediate consequence of that stability.

A special case of this theorem which is of independent interest is when the graph is the complete graph $K_n$.
In this case, the voter model is equivalent to the Moran model of population dynamics.
Our result has the interpretation that the eventual dominant allell is stable with respect to the initial population types, but sensitive to noise affecting the reproduction dynamics of the model.

\paragraph{Generalizations.}
The theorem and the proof extend with no significant changes to the case of weighted graphs (electrical networks).
Here each undirected edge has a weight $w_e$.
For each neighbouring pair $e=(x,y)$, the vertex $x$ adopts the opinion of $y$ at rate $w_e$, and $y$ adopts the opinion of $x$ at the same rate.
However, the rates in the two directions must be equal for our proof to carry through.
This can be seen to be necessary, since for general asymmetric rates the result may fail.
Even in the natural case where each vertex is activated at rate 1, and copies the opinion of a uniform neighbour the consensus might not be sensitive to noisy dynamics.
For example, in this model on the star graph, the consensus is strongly correlated to the initial opinion of the centre vertex, and this correlation is maintained when noise is added to the dynamics.

\paragraph{Proof overview.}

  The proof of Theorem~\ref{t:consensus_opinion} relies on the classical duality relation between the voter model and coalescing random walks.
  The same clock processes can be used to define coalescing random walks on $G_n$.
  If one fixes a realization of the clock processes $\sP$ and considers coalescing random walks going backwards in time and moving in the opposite direction from the activated edges, the consensus opinion coincides with the opinion at time zero from the last walker remaining.
  This observation alone implies noise stability with respect to the initial condition.

  Via this relation with coalescing random walks, one can also obtain that the position of the site that realizes the consensus opinion on $G_{n}$ is distributed according to the invariant measure of the random walk on $G_{n}$.
  This hints at the proof of the second statement about noise sensitivity with respect to clocks.
  Here one needs to study the joint evolution of two highly dependent random walks in $G_{n}$, using clock ticks from $\sP$ and $\sP^{\epsilon}$, respectively.
  One considers these walkers as evolving on $G_{n} \times G_{n}$ in a suitably constructed way.
The proof is then concluded by describing the invariant measure of this coupled dynamics and verifying that the probability that both random walks end in the same position vanishes in the limit as $n$ increases.

\paragraph{Related works.}

The notion of noise sensitivity of boolean functions was formalized by Benjamini, Kalai and Schramm \cite{bks}, following earlier related works of Hastad \cite{Hastad}.
It has developed into a central tool in the study of boolean functions, see e.g.\ O'Donnell's book \cite{ODonnell}.
Noise sensitivity and stability have proved to be especially useful in percolation theory, where the natural inputs are independent variables, and many of the previous  applications of that theory had been in the context of percolation.
In particular, noise sensitivity has proved crucial to our understanding of dynamical percolation.
See e.g. Garban and Steif's book \cite{GSbook} and references therein.

The consensus opinion of the voter model is known to have the same distribution as the independent initial opinions.
Obtaining bounds on the time when consensus is reached is more delicate, and also relies on the duality relation between voter model and coalescing random walks.
The consensus time of the voter model is bounded by the coalescence time of coalescing random walks.
Oliveira~\cite{oliveira} proves that, for graphs with constant degree, this expected coalescing time can be bounded up to a universal constant by the maximum of the hitting time in the graph over all starting vertices and targets.
This implies that the consensus time in the voter model is of order at most $n^{2}$.

\paragraph{Acknowledgements.}
We thank Itai Benjamini for raising the questions studied in this paper.
GA counted on the support of the Israel Science Foundation through Grant 957/20.
OA is supported in part by NSERC.
RB thanks the Mathematical Institute of Leiden University for support.
RP was supported in part by the Israel Science Foundation Grant 2566/20.

\section{Preliminaries}
\label{sec:prelim}

This section contains a short review of the voter model and its connection to coalescing random walks.
We highlight the main facts that will be used throughout the paper.

\paragraph{The voter model.}
Let $G_{n}=(V_n, E_{n})$ be a sequence of connected finite graphs.
The voter model with opinions in $S$ is a Markov chain on $S^{V_n}$, where a voter changes their opinion to $x\in S$ at rate equal to the number of neighbours with opinion $x$.
The voter model on $G_{n}$ can be obtained via a graphical construction that we define now.
Let $\vec{E}_n$ be the set of oriented edges of $G_n$, containing each edge in both directions.
Let $\sP$ denote a Poisson point process on $\vec{E}_{n} \times \R_{+}$ with intensity $\mu \otimes \lambda$, where $\mu$ is the counting measure on $\vec{E}_{n}$, and $\lambda$ denotes the Lebesgue measure in $\R_{+}$.
Given an initial configuration $\eta_{0} \in S^{V_n}$, the evolution proceeds by copying opinions at the time of the clock ticks: if the oriented edge $(u,v)$ has a clock tick at time $t$, then site $v$ copies the opinion of $u$.
It is easy to see that this almost surely defines a unique process $(\eta_{t})_{t \geq 0}$.

The voter model has absorbing states: the constant configurations with any opinion.
Furthermore, on any connected finite graph the model almost surely reaches one of these configurations.
An easy way to see this is by noting that the number of vertices with any opinion is a bounded martingale, and so must converge.

The voter model with general opinion set $S$ can be understood in terms of the voter model with Boolean opinions by noting that for any opinion $a\in S$, the set of vertices with opinion $a$ evolves as a voter model.
(These are not independent for different $a$'s but this is sufficient to deduce our results for a general opinion set from the result for Boolean opinions.)

\paragraph{Coalescing random walks.}
In this model, each site initially receives an independent random walk that jumps to each neighbour with rate one.
When any two of these walkers meet, they coalesce and move together from that time onward.

One way to define this model precisely is using a graphical construction.
Given Poisson processes $\sP$ as before, a random walk from $x$ is constructed by jumping from $u$ to $v$ if the walk is at $u$ and the clock for directed edge $(u,v)$ ticks.
It is clear that each walker performs a random walk, and that once two walkers are at the same vertex they stay together from that time on.
Since any two walkers on a finite connected graph will meet at some finite time, there is some time at which all walkers have coalesced.
This is called the coalescence time, which we denote by $T_C$.



\paragraph{Duality.}
The relation between the voter model and coalescing random walks is an easy consequence of the graphical constructions.
It makes it possible to use coalescing random walks to determine the consensus opinion via realizations of the graphical construction $(\eta_{0}, \sP)$.
Given $t \geq 0$, let us denote by $\big(X^{t}_{i}(s)\big)_{s \leq t}$ the random walk that is at position $i$ at time $t$ and runs backwards in time up to time zero by using the clock ticks from $\sP$: if an oriented edge $(u,v)$ has a clock tick and the random walk is on $v$, then it jumps to site $u$.
For each $t$, the walks running backward in time are coalescing random walks.

To see the relation with the voter model, one should trace the origin of the opinion of $x$ at time $t$.
Going backward in time, the opinion jumps from $v$ to $u$ if the edge $(u,v)$ has a clock tick.
Thus the origin of the opinions form coalescing random walks when viewed backward in time.

It is clear that the coalescence time of the backward random walks has the same distribution as that for walks moving forward in time.
Since $T_C$ is almost surely finite, the probability that the coalescence time for the walks $\big(X^{t}_{i} \big)_{i \in [n]}$ is at most $t$ is the same as $\P(T_C\le t)$ which tends to $1$ as $t\to\infty$.
Now, if $t$ is large enough such that the coalescing random walks with the collection $\big(X^{t}_{i} \big)_{i \in [n]}$ have coalesced by time $t$,
then for all $i$, we have $X_i^t(0)$ is the same vertex, and
\begin{equation}\label{eq:consensus}
  f_{n}(\eta_{0}, \sP) = \eta_{0}(X_{i}^{t}(0)).
\end{equation}
(Recall $f_{n}(\eta_{0}, \sP)$ denotes the consensus opinion of the voter model.)
Since $T_C$ is finite almost surely, the equation above also implies that
\begin{equation}
  f_{n}(\eta_{0}, \sP) = \lim_{t \to \infty} \eta_{0}(X_{i}^{t}(0)),
\end{equation}
which in turn provides the equality in distribution
\begin{equation}\label{eq:equality_distribution}
  f_{n}(\eta_{0}, \sP) \sim \eta_{0}(X),
\end{equation}
where $X$ is distributed according to the invariant measure of the random walk $X_{1}$.
If the initial opinions are independent, then $f_n$ has the same law as the initial opinions.

\paragraph{Coalescing random walks and noise.}

When studying noise sensitivity with respect to the clock rings, one needs to consider a pair of coupled random walks that evolve according to the pair of graphical constructions $(\sP, \sP^{\epsilon})$.
These random walks are coupled as follows:
If they are at distinct vertices they jump independently at rate 1 across each edge;
If they are at the same vertex, they jump together across each edge at rate $1-\epsilon$, and each of the walks jumps alone across the edge at rate $\epsilon$, leaving the other walk unchanged.

Proceeding in an analogous way as in the deduction of~\eqref{eq:equality_distribution}, we obtain the equality of the joint distribution
\begin{equation}\label{eq:equality_distribution_2}
\big( f_{n}(\eta_{0}, \sP), f_{n}(\eta_{0}, \sP^{\epsilon}) \big) \sim \big( \eta_{0}(X), \eta_{0}(X^{\epsilon}) \big),
\end{equation}
where $(X,X^{\epsilon})$ is distributed according to the invariant measure of a pair of coupled random walks evolving with the Poisson processes $(\sP, \sP^{\epsilon})$.

\section{Proof of Theorem~\ref{t:consensus_opinion}}

We now proceed with the proof of Theorem~\ref{t:consensus_opinion}.
We first verify noise stability with respect to the initial condition and then consider noise sensitivity with respect to clock ticks.

\paragraph{Noise stability for the initial condition.}
This is the easier part of the result.
It follows directly from \eqref{eq:consensus} by noticing that
\begin{equation}
  \P_{p} \Big[f_{n}(\eta_{0}, \sP) \neq f_{n}(\eta_{0}^{\epsilon}, \sP) \Big] = \P_{p} \Big[ \lim_{t \to \infty} \eta_{0}(X_{1}^{t}(0)) \neq \lim_{t \to \infty}\eta_{0}^{\epsilon}(X_{1}^{t}(0)) \Big]
  \leq \epsilon,
\end{equation}
since the two are equal if $\eta_0(X^t_1(0))$ has not been resampled.

\paragraph{Noise sensitivity for the clocks.}
  We now prove the statement regarding noise sensitivity with respect to the clock ticks in the graphical construction.
  We use the relation established in~\eqref{eq:equality_distribution_2} to obtain an expression for the distribution of the vertices that realize the consensus opinion in the coupled graphical constructions $(\sP, \sP^{\epsilon})$.

  Consider a pair $(X_{t}, X_{t}^{\epsilon})$ of random walks that evolve according to the pair of Poisson clocks $(\sP, \sP^{\epsilon})$.
  From the construction of the coupled pair of Poisson clocks, we see that $(X_{t}, X_{t}^{\epsilon})$ evolves on $G_{n} \times G_{n}$ with rates given by the $n^2\times n^2$ matrix
\begin{equation}
  Q\big((u,v), (u',v') \big) = \begin{cases}
    1-\epsilon, \text{ if } u=v, u'=v', \text{ and } u \sim u'; \\
    \epsilon, \text{ if } u=v=v' \text{ and } u \sim u'; \\
    \epsilon, \text{ if } u=u'=v \text{ and } v \sim v'; \\
    1, \text{ if } u \neq v, u=u', \text{ and } v \sim v'; \\
    1, \text{ if } u \neq v, v=v', \text{ and } u \sim u'; \\
    0, \text{ otherwise}.
  \end{cases}
\end{equation}

The key observation is that we can write explicitly the invariant measure of the random walk on $G_{n} \times G_{n}$ with transition rates given by the matrix $Q$.
Define the function $\nu$ (viewed as a $1\times n^2$ vector) on $G_{n} \times G_{n}$ by
\begin{equation}
  \nu(u,v) = \textbf{1}_{u = v}+ \epsilon \textbf{1}_{u \neq v}.
\end{equation}
We claim that $Q$ is reversible w.r.t. $\nu$.
Indeed, consider the quantity \newline $\nu(u,v) Q((u,v),(u',v'))$.
\begin{itemize}
\item If $u=v$ and $u'=v'$ this is $1\cdot(1-\epsilon)$ if $u\sim u'$ and 0 otherwise.
\item If $u=v$ and $u'\neq v'$, this is $1\cdot\epsilon$ if one of $\{u',v'\}$ is $u$ and the other is a neighbouring vertex (otherwise 0).
\item If $u\neq v$ and $u'=v'$, this is $\epsilon\cdot 1$ if one of $\{u,v\}$ is $u'$ and the other is a neighbouring vertex (otherwise 0).
\item If $u\neq v$ and $u'\neq v'$, this is $\epsilon\cdot 1$ as long as one of the pairs $(u,u')$, $(v,v')$ is equal and the other are neighbours.
\end{itemize}
In all cases, this is the same if $(u,v)$ and $(u',v')$ are swapped.


This implies that the invariant measure $\pi$ for the process $(X_{t}, X_{t}^{\epsilon})_{t \geq 0}$ is proportional to $\nu$, i.e. $\pi = \nu/(n+(n^2-n)\epsilon)$.
In particular $\pi(X=X') = \frac{1}{1+(n-1)\epsilon}$.

Suppose $\eta_0$ has Bernoulli$(p)$ entries.
We now calculate:
\begin{equation}
  \begin{split}
    \E_{p}\Big( f_{n}(\eta_{0}, \sP) f_{n}(\eta_{0}, \sP^{\epsilon}) \Big)
    & = \E_\pi \E_{p}\Big( \eta_{0}(X), \eta_{0}(X^{\epsilon}) \Big) \\
    & = p^2 \pi(X\neq X') + p \pi(X=X') \\
    & = p^2 + \frac{p-p^2}{1+(n-1)\epsilon}.
  \end{split}
\end{equation}
Thus $f_{n}(\eta_{0}, \sP)$ and $f_{n}(\eta_{0}, \sP^{\epsilon})$ are asymptotically independent if and only if $\epsilon n\to\infty$, as claimed.

\section{A problem}

Our theorem states that the consensus is sensitive to noisy dynamics.
Moreover, this sensitivity does not depend on the graph on which the process lives.
Consider a new model, which we call the \textbf{majority-voter} model.
Here, the voter dynamics are run for some pre-determined time $T$.
At that time, an election occurs, and the outcome is the majority of the opinions at time $T$.
This interpolates between simple majority for $T=0$ and the consensus in the voter model for $T\to\infty$.
A natural question is how large does $T$ needs to be so that the majority-voter model is sensitive to noisy dynamics.
It seems that the result would depend on the underlying graph.

\bibliographystyle{plain}
\bibliography{mybib}

\begin{thebibliography}{1}

\bibitem{bks}
Itai Benjamini, Gil Kalai, and Oded Schramm.
\newblock Noise sensitivity of {B}oolean functions and applications to
  percolation.
\newblock {\em Publications Math{\'e}matiques de l'Institut des Hautes
  {\'E}tudes Scientifiques}, 90(1):5--43, 1999.

\bibitem{GSbook}
Christophe Garban and Jeffrey~E Steif.
\newblock Noise sensitivity and percolation.
\newblock {\em Probability and statistical physics in two and more dimensions},
  15:49--154, 2012.

\bibitem{Hastad}
Johan H{\aa}stad.
\newblock Some optimal inapproximability results.
\newblock {\em Journal of the ACM (JACM)}, 48(4):798--859, 2001.

\bibitem{liggett}
Thomas~M Liggett.
\newblock Interacting particle systems - an introduction.
\newblock In {\em School and Conference on Probability Theory}, 2004.

\bibitem{ODonnell}
Ryan O'Donnell.
\newblock {\em Analysis of {B}oolean functions}.
\newblock Cambridge University Press, 2014.

\bibitem{oliveira}
Roberto Oliveira.
\newblock On the coalescence time of reversible random walks.
\newblock {\em Transactions of the American Mathematical Society},
  364(4):2109--2128, 2012.

\end{thebibliography}

\end{document}